\documentclass[%
 reprint,
showpacs,
 amsmath,amssymb,
 aps,
]{revtex4-1}

\usepackage{graphicx}
\usepackage{dcolumn}
\usepackage{bm}

\allowdisplaybreaks

\usepackage{mathrsfs}

\usepackage[usenames,dvipsnames,svgnames,table]{xcolor}
\colorlet{linkequation}{blue}
\usepackage[colorlinks,urlcolor=lightgray,citecolor=black,pdfborder={0 0 0}]{hyperref}

\renewcommand{\i}{\mathrm{i}}
\newcommand{\re}{{\text{re}}}
\newcommand{\im}{{\text{im}}}

\begin{document}

\title{An algebraic approach to the Kuramoto model}

\author{Lyle Muller}
\author{J\'an Min\'a{\v c}}
\author{Tung T. Nguyen}
\affiliation{Department of Mathematics, Western University \\ London, Ontario, Canada N6A 3K7}

\date{\today}


\begin{abstract}
We study the Kuramoto model with attractive sine coupling. We introduce a complex-valued matrix formulation whose argument coincides with the original Kuramoto dynamics. We derive an exact solution for the complex-valued model, which permits analytical insight into individual realizations of the Kuramoto model. The existence of a complex-valued form of the Kuramoto model provides a key demonstration that, in some cases, re-formulations of nonlinear dynamics in higher-order number fields may provide tractable analytical approaches. 
\end{abstract}


\maketitle


The dynamics of networks with many nodes and connections poses difficulties for mathematical treatment. While linear dynamics on a network of interacting units is straightforward to handle with matrix algebra, once nonlinearity is introduced at individual units, the dynamics of the system becomes difficult to study analytically. Here, we consider the Kuramoto model (KM), a paradigmatic example of nonlinear network dynamics describing the synchronization of coupled oscillators \cite{Acebron05}. This model provides a canonical description of synchronization in nature, where populations of interacting units (from neurons to Josephson junctions and fireflies \cite{Strogatz03}) coordinate the timing of their behavior in the absence of a central coordinator \cite{Strogatz01}. Because of its wide applicability throughout many systems, the KM is a central consideration in the study of nonlinear dynamics.

Analytical approaches to the KM started with the original work of Kuramoto, who introduced the model in \cite{Kuramoto75}. By changing to a rotating coordinate frame and passing to a continuum description, Kuramoto provided a description of the population dynamics in the infinite limit and the transition to synchrony at a critical coupling strength. Since this pivotal work, much research has gone on to analyze the dynamics of KM. In recent years, several studies have focused on introducing complex versions of the KM \cite{Roberts08,vanMieghem09}, with the goal of providing some analytical insight; however, finding a complex KM that permits an exact solution for the evolution of this system has proven evasive.

In this Letter, we provide a complex-valued matrix formulation of the KM whose argument corresponds to the original KM. We derive an explicit analytical solution for the dynamics using the eigenspectrum of the adjacency matrix. We then use this solution to study synchronization in individual realizations of the model. The existence of a complex-valued version of the KM permitting an analytical solution provides a key demonstration for potential analytical approaches to complex nonlinear dynamics at finite scales. The complex-valued KM and its exact solution provide an example that, in some cases, re-formulations of nonlinear dynamics in different number fields may provide opportunities for analytical descriptions of nonlinear systems.

We start with the standard Kuramoto model (KM) on a general network of $N$ nodes:

\begin{equation}\label{KM}
\dot{\theta}_i = \omega_i + \kappa \sum_{j=1}^{N} a_{ij} \sin( \theta_j - \theta_i )\,,
\end{equation}

\noindent where $\theta_i(t) \in [-\pi,\pi]$ is the state of oscillator $i \in [1,N]$ at time $t$, $\omega_i$ is the intrinsic angular frequency, $\kappa$ scales the coupling strength, and element $a_{ij} \in \{0,1\}$ represents the connection between oscillators $i$ and $j$. The standard sine coupling in the interaction term causes phases of two connected oscillators $i$ and $j$ to attract, to an extent depending on the homogeneous coupling strength $\kappa$. Here, we first consider the KM defined on an undirected ring graph $\mathfrak{G}_{RG}$, where nodes are arranged on a one-dimensional ring with periodic boundary conditions and connected to their $k$ neighbors in each direction (later considering Erd\H{o}s-R\'enyi and Watts-Strogatz graphs). We also consider first the case of homogeneous intrinsic frequencies $\omega_i = \omega\,\forall\,i \in [1,N]$ (with inhomogeneous case in Supplement \cite{Supplement}). For $k = \lfloor N/2 \rfloor$, $\mathfrak{G}_{RG}$ corresponds to the fully connected graph on $N$ nodes, $K_N$. It is important to note that the approach described here is general to the specific graph model we consider.

\begin{figure}[t]
\includegraphics[width=0.98\columnwidth]{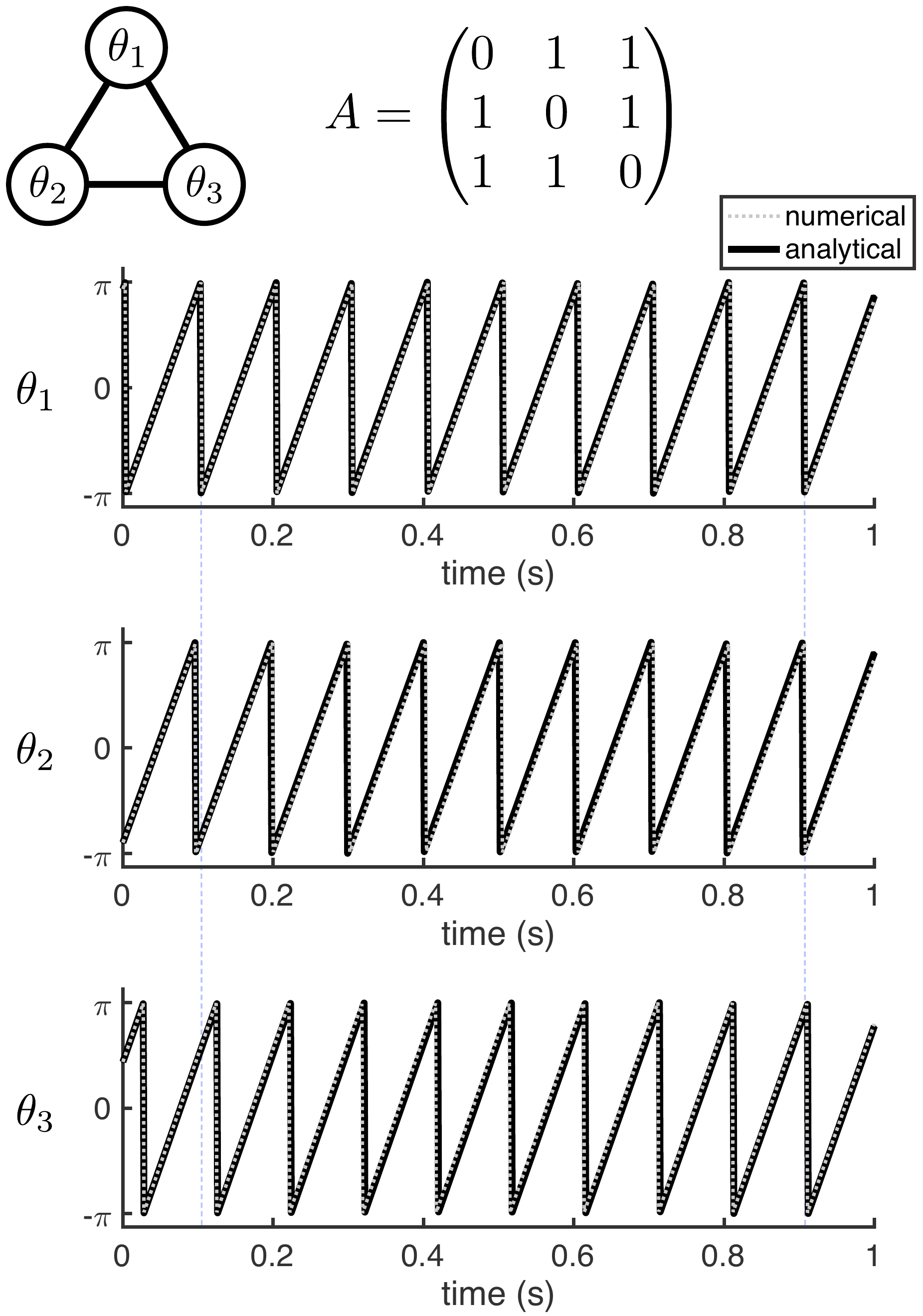}
\caption{\label{Fig_Coupling} 
Synchronization dynamics in a three-node Kuramoto model (KM). In this example, $N = 3$, $\kappa = 1$, $\omega/2\pi = 10\,\text{Hz}$, and the network is fully connected (top left), so the adjacency matrix $\bm{A}$ is the complete graph $K_3$ (top middle). Compared are the numerical result (dotted gray line) and the analytical expression (solid black line) for the state variables $\theta_1$ (top), $\theta_2$ (middle), and $\theta_3$ (bottom). Plotted is one second of simulation ($dt = 0.001\,\text{s}$). The difference between the numerical result and analytical prediction was verified to remain less than $\pi/16$ in simulations as long as 10 seconds.}
\end{figure}

Starting with Equation \eqref{KM}, we can change to a rotating coordinate frame \cite{Strogatz88} and set $\omega = 0$ without loss of generality. We then subtract an additional imaginary component in the interaction term:

\begin{equation}\label{KMmod}
\dot{\theta}_i = \gamma \sum_{j=1}^{N} a_{ij} \big[ \sin( \theta_j - \theta_i ) - \i \cos( \theta_j - \theta_i ) \big]\,.
\end{equation}

\noindent We note this expression now implies $\theta_i \in \mathbb{C}$ and requires a scaling of the coupling strength ($\gamma = 2\kappa/\pi$); as we will show, the resulting complex-valued system admits a solution whose argument agrees with the original KM.

Multiplying both sides of Equation \eqref{KMmod} by $\i$, we have

\begin{equation} 
\i \dot{\theta}_i = \gamma \sum_{j=1}^{N} a_{ij} \big[ \i \sin( \theta_j - \theta_i ) + \cos( \theta_j - \theta_i ) \big]. 
\end{equation}

\noindent Applying Euler's formula, we can rewrite the system as:

\begin{equation}
\i \dot{\theta}_i = \gamma \sum_{j=1}^{N} a_{ij} e^{\i (\theta_j-\theta_i)}=\gamma e^{-\i\theta_i} \sum_{j=1}^{N} a_{ij} e^{\i\theta_j}.
\end{equation}

\noindent This equation now results in the following matrix form:

\begin{equation}
\bm{\dot{\theta}} = \frac{\gamma}{\i}\,\text{diag}[ e^{-\i\bm{\theta}} ]\,\bm{A} e^{\i\bm{\theta}}\,,
\end{equation}

\noindent where $\bm{\theta} \in \mathbb{R}^N$ is the state vector across nodes, and $\bm{A} \in \{0,1\}^{N \times N}$ is the adjacency matrix encoding the connections between nodes in the network. We can then utilize the fact that the inverse of the matrix $\text{diag}[ e^{-\i\bm{\theta}} ]$ is $\text{diag}[ e^{\i\bm{\theta}} ]$ to arrive at

\begin{equation}\label{eqDiag}
\text{diag}[ e^{\i\bm{\theta}} ] \bm{\dot{\theta}} = \frac{\gamma}{\i} \bm{A} e^{\i\bm{\theta}}\,.
\end{equation}

Now, using the fact that 

\begin{equation}
\text{diag}[ e^{\i\bm{\theta}} ] \bm{\dot{\theta}} = \frac{1}{\i} \frac{d}{dt} e^{\i\bm{\theta}}\,,
\end{equation}

\noindent mentioned in \cite{vanMieghem09}, we have that

\begin{equation}
\frac{d}{dt} e^{\i\bm{\theta}} = \gamma \bm{A} e^{\i\bm{\theta}}
\end{equation}

\noindent and letting $\bm{x} = e^{\i\bm{\theta}}$, we have

\begin{equation}
\bm{\dot{x}} = \gamma \bm{A} \bm{x}
\end{equation}

\noindent whose general solution is

\begin{equation}\label{eqSol}
\bm{x}(t) = e^{\gamma t \bm{A}} \bm{x}(0)\,.
\end{equation}

\noindent We can now utilize a diagonalization of the adjacency matrix, $\bm{A} = VDV^{-1}$, to rewrite the solution as

\begin{equation}
\bm{x}(t) = \bm{V} e^{\gamma t \bm{D}} \bm{V^{-1}} \bm{x}(0)\,,
\end{equation}

\noindent and the problem of computing the matrix exponential is now reduced to computing the standard exponential function on the diagonal entries of $\bm{D}$. We note that the adjacency matrix $\bm{A}$ is always diagonalizable for the undirected graphs considered here.

Finally, as noted above, Equation \eqref{eqSol} can be directly related to the original KM. Let $\bm{\theta}=\bm{\theta}_{\re}+\i \bm{\theta}_{\im}$ be the decomposition of $\bm{\theta}$ into the real and imaginary parts. Then we have

\begin{equation}
\bm{x}=e^{\i \bm{\theta}_{\re}- \bm{\theta}_{\im}}=e^{-\bm{\theta_{\im}}} e^{\i \bm{\theta}_{\re}}
\end{equation}

\noindent $\bm{\theta_{\re}}$ is thus the argument of the analytical solution $\bm{x}$. In particular, we can take $\bm{\theta_{\re}} \in [-\pi, \pi]$. We will use this solution to compare with the numerical integration of Equation \eqref{KM}.

\begin{figure}[t]
\includegraphics[width=\columnwidth]{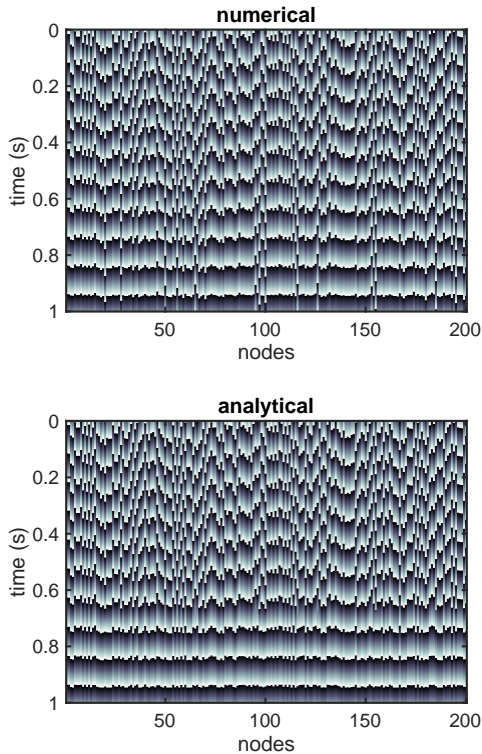}
\caption{Numerical and analytical realizations of the KM. Plotted are the state variables ($\theta_i\,\forall\,i \in [1,N]$) produced by the model ($N = 200$, $k = 100$, resulting in the fully connected complete graph $K_N$), with different nodes on the horizontal axis and time on the vertical. Dark colors indicate values close to $-\pi$, and light colors indicate values close to $\pi$.}
\end{figure}

What remains in deriving an expression for the resulting temporal dynamics is to study the eigenspectrum of $\bm{A}$. To do this, we note the adjacency matrix of the ring graph $\mathfrak{G}_{RG}$ is by definition a circulant matrix \cite{Davis79}, and we can use the circulant diagonalization theorem (CDT) to obtain its eigenspectrum analytically. The CDT states that all circulants $c_{ij} = \text{circ}(c_j)$, where $\text{circ}(c_j)$ is the circulant matrix constructed from the generating vector $c_j$, are diagonalized by the same unitary matrix $\bm{U}$ with components

\begin{equation}
u_{rs} = \frac{1}{\sqrt{N}} \exp\left[ - \frac{2\pi \i}{N} (r-1)(s-1) \right] ,
\end{equation}

\noindent $r,s \in [1,N]$ and that the $N$ eigenvalues are given by

\begin{equation}
E_r(\mathbf{C}) = \sum\limits_{j=1}^{N} c_j \exp\left[ - \frac{2\pi \i}{N} (r-1)(j-1) \right]\,.
\end{equation}

\noindent Using these expressions, we can then evaluate Equation \eqref{eqSol} in terms of the eigenspectrum of $\bm{A}$ to obtain a fully analytical evaluation of $\bm{x}(t)$, which we can then compare to numerical integration of the original KM.

We can now compare the argument of the analytical expression Equation \eqref{eqSol} with the result obtained by numerical integration of Equation \eqref{KM}. As a first example, we considered a small network of $N = 3$ nodes and $k = 1$, resulting in the complete graph $K_3$. Figure 1 (top left) shows the considered network along with the state variables $\theta_i$ at each node. Initial conditions $\theta_i(0)$ were selected randomly with uniform distribution~~$\mathcal{U}(-\pi,\pi)$. Starting from these random initial conditions and a homogeneous coupling strength $\kappa = 1$, this small network synchronizes over the course of the one second simulation (cf.\,dotted blue vertical lines across bottom three panels, Fig.~1). Equations were integrated using a forward Euler method with fine temporal precision ($dt = 0.001\,\text{s}$) and were compared with results from a Runge-Kutta method \cite{Dormand80}. Additional simulations were run at very high temporal precision (up to $dt = 10^{-6}\,\text{s}$) to ensure the numerical validity of these results. Simulation code to reproduce all figures in this work will be made availabe at our GitHub site (\href{http://mullerlab.github.io}{http://mullerlab.github.io}).

We next compare the analytical and numerical results for larger networks of $N = 200$ nodes on a ring graph with $k = 100$. To create a synchronization transition during a one second simulation interval, we scaled the coupling strengths accordingly ($\kappa = 6/N$). Figure 2 shows the results of this comparison between the numerical (top) and analytical (bottom) evaluations. The numerical and analytical evaluations exhibit similar dynamics, from the macroscopic synchronization to the specific trajectories of individual oscillators. Interestingly, a small fraction of the nodes in the numerical simulation remain counterphase to the rest of the population, while this behavior does not occur in the analytical version.

\begin{figure}[b]
\includegraphics[width=\columnwidth]{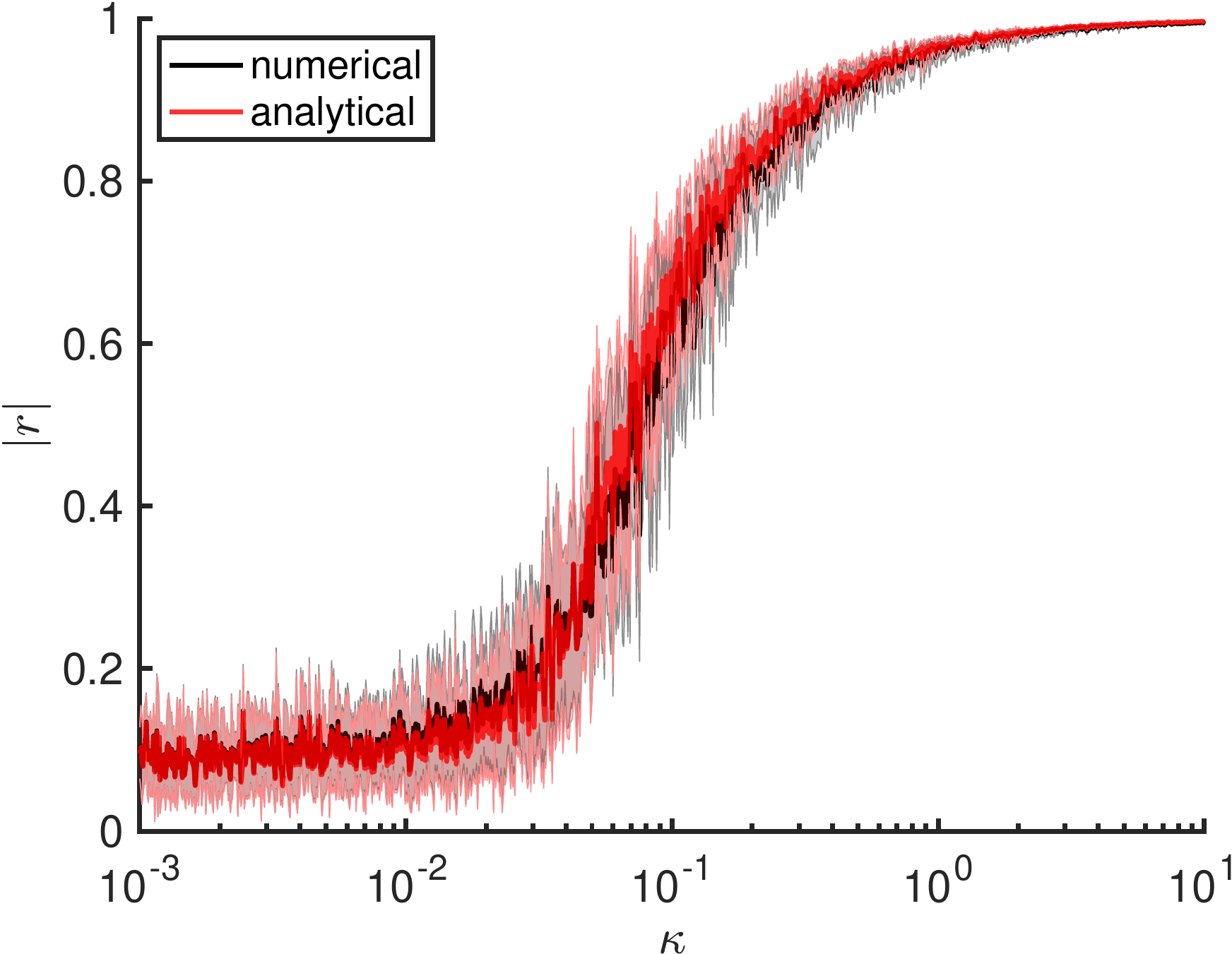}
\caption{Synchronization in numerical and analytical realizations of the KM ($N = 200$, $k = 100$). The modulus of the synchronziation order parameter $r$ is plotted as a function of the coupling strength $\kappa$ for the numerical (black line) and analytical (red line) evaluations of the model. Solid lines represent the average $|r|$ over 10 1-second simulations (with the same random initial conditions for the numerical and analytical evaluation in each realization) at a value of $\kappa$ on a logarithmic space of 1,000 points between $10^{-3}$ and $10^{1}$. Shaded areas represent standard deviation over realizations.}
\end{figure}

To understand this point further, we systematically studied synchronization in the KM. To quantify the extent of synchronization, we use a standard approach to study the sum over oscillators:

\begin{equation}
r(t) = \frac{1}{N} \sum_{j=1}^{N} e^{\i\theta_j(t)}
\end{equation}

\noindent where $r(t) \in \mathbb{C}$, $\lvert r(t) \rvert \in [0,1]$ is the order parameter at time $t$, and $\theta_j(t) \in [1,N]$ is the phase of oscillator $j$ at time $t$ (n.b.~$\theta_{j,\re}(t)$ in the case of the analytical expression Equation \ref{eqSol}). Figure 3 shows the time-average order parameter over 1-second simulations. As $\kappa$ increases in the KM, the order parameter begins at a low value (desynchronized state) and increases until approaching unity (synchronized state) at values of $\kappa$ ranging from 1 to 10. As observed in Figure 2, the numerical and analytical versions of the KM exhibit very similar macroscopic synchronization dynamics.

Lastly, we considered the numerical and analytical solutions of the KM on undirected Erd\H{o}s-R\'enyi and Watts-Strogatz random graphs. For each realization of a random graph, we obtain a numerical estimate of the eigenspectrum of its adjacency matrix $\bm{A}$ to use in the analytical expression Equation \eqref{eqSol}. We first considered the KM on an Erd\H{o}s-R\'enyi random graph ($\mathfrak{G}_{ER}$), which displays synchronization dynamics similar to those previously observed (Fig.~4, top). We then considered the KM on a Watts-Strogatz network ($\mathfrak{G}_{WS}$), which is defined as a ring graph $\mathfrak{G}_{RG}$ where each node is first connected to its $k$ neighbors in each direction and each edge is rewired to another node with uniform probability $q$ \cite{WattsStrogatz98}. The KM on $\mathfrak{G}_{WS}$ displays non-trivial spatiotemporal dynamics before converging to the synchronized state (Fig.~4, bottom middle). Importantly, these spatiotemporal dynamics were also well described by the analytical expression introduced here (Fig.~4, bottom right).

In this work, we have introduced a complex-valued formulation of the KM whose argument corresponds to the original Kuramoto dynamics. This formulation permits an exact analytical solution for individual realizations of the KM. Here, we have first considered the case of homogeneous intrinsic frequency; however, this approach generalizes to the inhomogeneous case (Supplement \cite{Supplement}). We then compared the analytical version to numerical integration of the KM equations and found the analytical version displays similar dynamics, from the macroscopic synchronization behavior to the specific trajectories of individual oscillators.

Previous research, including an inspiring technical report written by van Mieghem \cite{vanMieghem09}, has studied complex-valued formulations of the KM \cite{Roberts08,Cadieu10}; however, no exact analytical expression was previously obtained. In particular, the expression derived in \cite{vanMieghem09} was noted to hold only for the repulsive cosine variant of the KM, and the linear reformulation in \cite{Roberts08} requires tuning a parameter to create the correspondence to the original KM. Thus, the results reported in this Letter, while representing only an initial study utilizing the analytical expression in Equation \eqref{eqSol}, to the best of our knowledge represent the first analytical version of the Kuramoto model.

\begin{figure}[t]
\includegraphics[width=\columnwidth]{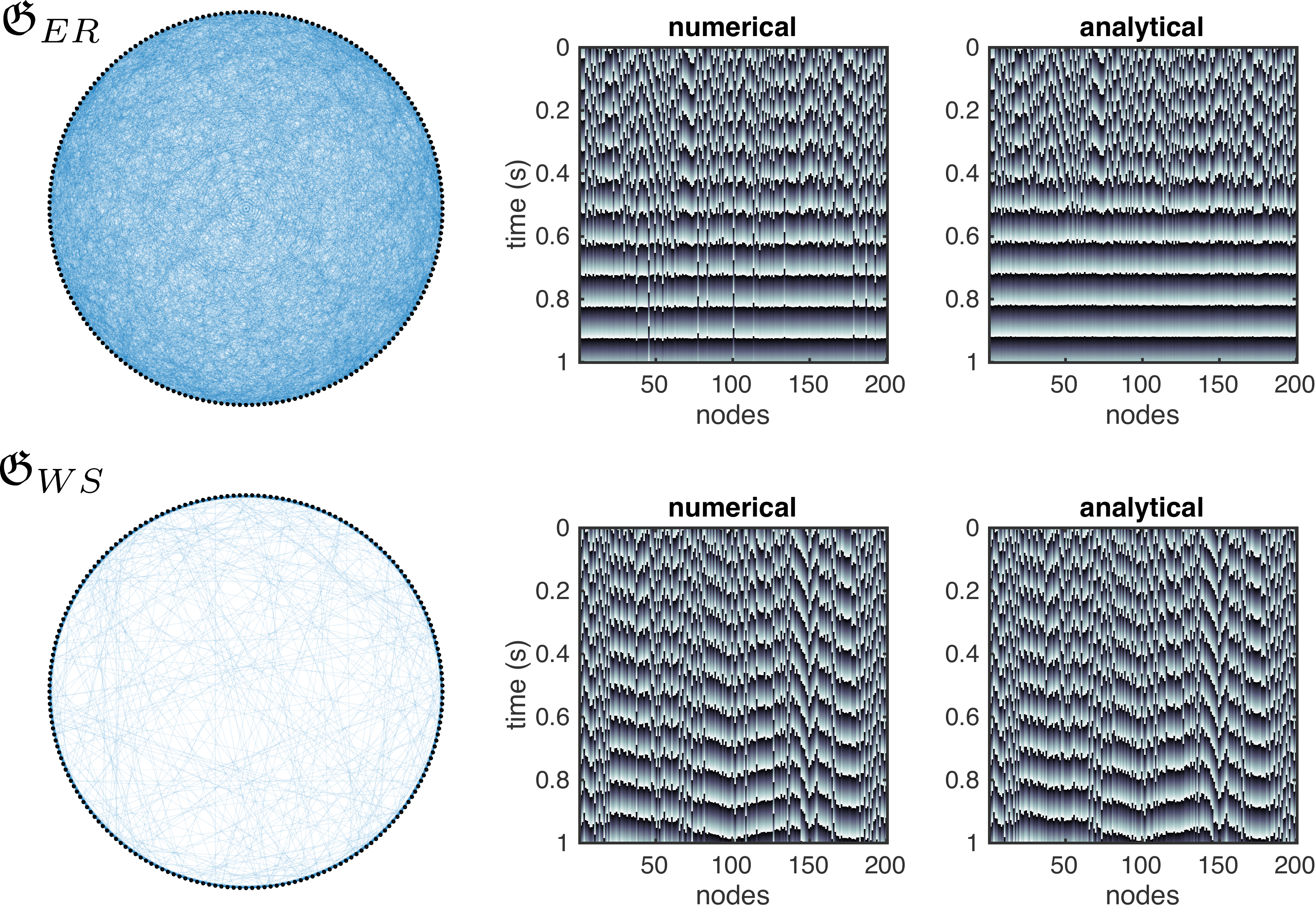}
\caption{Numerical and analytical realizations on the Erd\H{o}s-R\'enyi and Watts-Strogatz random graphs. (top left) The Erd\H{o}s-R\'enyi random graph $\mathfrak{G}_{ER}$ is plotted with nodes (black dots) and edges (blue lines) ($N = 200, p = 0.2$). (top right) The numerical integration of the KM (Eq.~\ref{KM}) and the analytical evaluation of Eq.~\eqref{eqSol} are plotted ($\kappa = 50/N$). As previously, dark colors indicate values close to $-\pi$, and light colors indicate values close to $\pi$. (bottom left) The Watts-Strogatz (WS) graph $\mathfrak{G}_{WS}$ ($N = 200$, $k = 10$, and rewiring probability $q = 0.1$). (bottom right) Numerical and analytical evaluations of the KM are plotted as above.}
\end{figure}

Importantly, we emphasize that the analytical version introduced here is valid at finite scales and for individual realizations of the KM. This analytical version allows future mathematical study of synchronization dynamics in networks with many nodes and connections, potentially using new tools from spectral graph theory, in addition to allowing one to obtain the future state of the system at an arbitrary future moment without numerical integration of the differential equations defined on the network. Because the KM has been extensively studied both as a model for neural dynamics \cite{Strogatz88,Mirollo90,Abbott93} and as a fundamental model for neural computation \cite{Kleinfeld01}, these results open up several possibilities for studying the connections between network structure, nonlinear dynamics, and computation. Recurrent connections have previously been shown to produce powerful computations through nonlinear interactions \cite{Sinha98}. The approach introduced here may have applications in understanding such recurrent interactions, which have been increasingly acknowledged to play an important and unexplained role in visual processing in the brain \cite{Kietzmann19}. Understanding more clearly the connection between networks and computation thus may have implications for fields such as neuroscience and beyond.

In this work, we have specifically studied the KM defined on a ring graph $\mathfrak{G}_{RG}$, whose highly regular structure permits analytical study of the eigenspectrum of its adjacency matrix. The regular structure of this graph means that its adjacency matrix belongs to the class of circulant matrices, whose eigenspectrum can be calculated analytically. Further, when $k$ is maximal, such that the number of neighbors to which a node is connected equals the rest of the graph, $\mathfrak{G}_{RG}$ corresponds to $K_N$, the complete graph on $N$ nodes. The KM defined on $K_N$, in turn, corresponds to the case of all-to-all connectivity first considered by Kuramoto \cite{Kuramoto75}. For these reasons, we chose to focus on the KM defined on $\mathfrak{G}_{RG}$ in this work. In previous work \cite{RudolphMuller14}, however, we have introduced an operator-based approach to the structure of random graphs. In future work, we aim to extend the present results to understand the connection between nonlinear dynamics and random graphs with various structural features.

The existence of a complex-valued formulation of the Kuramoto model that permits an analytical solution raises an important example in nonlinear dynamics. While Equation \eqref{KM}, which includes the sine coupling interaction term and the network adjacency matrix, appears analytically intractable, in this case a re-formulation in the complex domain provides an algebraic approach to the Kuramoto dynamics which serve as a canonical description for synchronization phenomena in nature. While this re-formulation may not generalize beyond the Kuramoto model, the ability of this approach to provide insight in this case suggests that representations in number fields beyond $\mathbb{R}$ may represent opportunities for future insight.


\begin{acknowledgments}
The authors wish to thank M Ly, E Bienenstock, and TP Coleman for insightful discussions. L.M. was supported by BrainsCAN at Western University through the Canada First Research Excellence Fund (CFREF) and by the NSF through a NeuroNex award (\#2015276). J.M. was supported by the Natural Sciences and Engineering Research Council of Canada (NSERC) grant R0370A01.
\end{acknowledgments}


\end{document}